\documentclass[a4paper,oneside,12pt]{amsart} 

\usepackage[T2A]{fontenc}
\usepackage[english]{babel}
\usepackage{amsmath, amssymb, amsfonts, amsthm}
\usepackage{enumitem}
\usepackage{graphicx}
\usepackage{tikz}
\usepackage{hyperref}
 \usepackage{dutchcal}
\usepackage{bbm}
\usepackage{geometry}

\newtheorem{lem}{Lemma}
\newtheorem{defin}{Definition}
\newtheorem{thm}{Theorem}
\newtheorem{rem}{Remark}

\def\XXint#1#2#3{{\setbox0=\hbox{$#1{#2#3}{\int}$ }
\vcenter{\hbox{$#2#3$ }}\kern-.6\wd0}}

\begin{document}

\title{Some remarks on $K-$ closedness of the couples of real Hardy spaces.}
\date{19 March 2014}

\author[Ioann Vasilyev]{Ioann Vasilyev}
\email{ioann.vasilyev@cyu.fr, milavas@mail.ru}
\thanks{Statements and Declarations: The author declares no competing interests.}

\begin{abstract}
In this paper $K$-closedness is proved in the case of the couple of real Hardy spaces  in the corresponding couple of Lebesgue spaces when both indexes are between $(n-1)/n$ and $\infty$. This means roughly that any measurable decomposition of an analytic function gives rise to an ``analytic'' decomposition with summands of roughly the same size. The proof uses Bourgain's method, the atomic decomposition for Hardy spaces and the subharmonic property of the gradient of a system of conjugate harmonic functions.
\end{abstract}

\maketitle

\section{Introduction}

It has been known for a long time that in the interpolation sense the scale of analytic Hardy classes on the unit circle acts in fact the same way as the scale $L^p$.  
If we consider the real method of interpolation, one sees that in that case everything is clear: the couple $(H^r,H^t)$ is $K-$closed in the couple $(L^r,L^t)$ for all $r,t \in (0,\infty].$ The last definition means that any decomposition of a function $f\in H^r+H^t$  can be ``transformed'' into  an analytic one with roughly the same values of summands. The most complicated cases here are when $r,t \in (0,1]\cup\{+\infty\}$. It was G. Pisier, who first proved that in the most common case. Consult \cite{bou} and \cite{svk} for the other methods. In the last paper one can also find the weighted case.

For $\infty>r>\frac{n-1}{n}$ the space $H^r(\mathbb R^n)$, as it is known, can be interpreted as the space of distributions $\phi$, such that $\phi , R_1 \phi ,...,R_n \phi \in L^r(\mathbb R^n)$, where $R_1,..R_n$ are the Riesz transformations. 
For here comes the embedding of $H^r(\mathbb R^n)$ into the space $L^r$, more precisely, into $(L^r (\mathbb R^n)\oplus...\oplus L^r(\mathbb R^n))$ (n+1 summand): 
$\phi \rightarrow (\phi, R_1 \phi ,...,R_n \phi)$.

 Interpolation properties of the scale $H^r(\mathbb R^n)$ are well understood, but concerning the question of $K$-closedness, the answer was known in only one case: $(H^1(\mathbb R^n),H^p(\mathbb R^n))$ is $K-$closed in $(L^1 (\mathbb R^n)\oplus...\oplus L^1(\mathbb R^n), L^p (\mathbb R^n)\oplus...\oplus L^p(\mathbb R^n))$. (For
$p>1$ Riesz transforms are bounded at $L^p$, so the case of both power greater than $1$, is also known, but is not interesting.)

In our paper we prove $K-$closedness in the case of the couples $(H^r(\mathbb R^n),L^p(\mathbb R^n))$ for $\frac{n-1}{n}<r<1<p<\infty$. 

\section{The main theorem}

Let us  first give a rigorous definition of $K-$ closedness.

\begin{defin}
The pair of (quasi) Banach spaces $(H_1,H_2)$ is $K$-closed in the pair $(L_1,L_2)$, if and only if, by definition, $H_i$ is closed in $L_i$, and for any function $f$ such that $f=f_1+f_2$, where $f_i\in H_i$, and for it's any decomposition $f=l_1+l_2$, where $l_i\in L_i$, there exists another decomposition $f=h_1+h_2$, where $h_i\in H_i$, and $\|h_i\|_{H_i}\leq C \|l_i\|_{L_i}$.
\end{defin}

We are now in position  to formulate the main result of our paper.

\begin{thm}
The couple of real Hardy spaces $(H^{p_1}(\mathbb R ^n),H^{p_2}(\mathbb R ^n))$ is $K-$closed in the couple of corresponding Lebesgue spaces $(L^{p_1}(\mathbb R ^n),L^{p_2}(\mathbb R ^n))$ for $\frac{n-1}{n} < p_1 < 1< p_2.$
\end{thm}

\begin{proof}

Let us first rewrite the statement of our theorem in more detail.

Let $f \in L^{p_1}(\mathbb R ^n), R_if=\alpha_i+\beta_i,$ where $R_i$ are the Riesz transforms and let $$\|(f, \alpha_1, \dots, \alpha_n)\|_{L^{p_1}\times \dots L^{p_1}}=\alpha,$$ $$\|(f, \beta_1, \dots, \beta_n)\|_{L^{p_2}\times \dots L^{p_2}}=\beta.$$ We also assume as we can that $\frac{n-1}{n}< p_1<1, p_2>1.$ 

So, we have to show that there exists functions $w\in H^{p_1}(\mathbb R ^n), v\in H^{p_2}(\mathbb R ^n)$ such that $f=w+v$ and, moreover, $\|w\|_{H^{p_1}}\leq C\alpha, \|v\|_{H^{p_2}}\leq C\beta.$

We shall start with the case $p_2=2$ and after that we shall prove the general case. Nevertheless, we remark that some parts of our proof in the case $p_2=2$ hold true for general $p_2$.

Consider the following set $$A:=\{x: \sup_{(y,t)\in \Gamma_x}|f*P_t(y,t)|>\lambda\},$$
where $\Gamma_x$ is the cone of the aperture $30n$ and with the vertex at the point $x$, and $$\lambda=\left(\frac{\beta^{p_2}}{\alpha^{p_1}}\right)^{\frac{1}{p_2-p_1}}.$$
Decompose the set $A$ with the help of the Whitney decomposition theorem
$$A=\bigcup_{i\ge 1}A_i.$$ 

Next, we estimate the Lebesgue measure $|A|$ of a set $A.$ 
\begin{lem}
The following estimate of $|A|$ holds
$$\left|A\right|^{\frac{p_1}{p_2}}\leq C\left(\frac{\beta^{p_1}}{\lambda^{p_1}} + \left(\frac{\alpha^{p_1}}{\lambda^{p_1}}\right)^{\frac{p_1}{p_2}} \right).$$
\end{lem}
\begin{proof}

Here is the main place of the proof, where we use the subharmonicity of the gradient of so-called system of conjugate harmonic functions (\cite {GCRF}) and the fact that $R_if=\alpha_i+\beta_i:$ (we introduce the following notation: $\frac{n-1}{n}=: \delta$)
$$\left|f* P_t(y,t)\right|^{\delta}\leq C\left(|\tilde\alpha|^{\delta}* P_t(y,t) + |\tilde\beta|^{\delta}* P_t(y,t) \right),$$
where $\tilde\alpha$ and $\tilde\beta$ are the functions, such that $\|\tilde\alpha\|_{L^{p_1}}\leq C \alpha, \|\tilde\beta\|_{L^{p_2}}\leq C \beta.$
Taking the supremum over a cone we get:
$$\sup_{(y,t)\in\Gamma_x}\left|f* P_t(y,t)\right|^{\delta}\leq C M |\tilde\alpha|^{\delta}+M |\tilde\beta|^{\delta},$$
where $M$ is the standard maximal function.
Raising the last line to the power $\frac{p_1}{\delta}$ we obtain the following inequality:
  $$\sup_{(y,t)\in\Gamma_x}\left|f* P_t(y,t)\right|^{p_1}\leq C(M |\tilde\alpha|^{\delta})^{\frac{p_1}{\delta}}+C(M |\tilde\beta|^{\delta})^{\frac{p_1}{\delta}}$$
We introduce one more auxiliary notation: $$A_N:=A\cap\{x\in \mathbb R^n: |x|\leq N\}.$$
We need that restriction because, generally speaking, $|A|$ can happen to be equal to $\infty$, while $|A_N|$ is finite. Next, we integrate out over the set $A_N$, using the Holder's inequality and the fact that $Mf$ a bounded operator acting from $L^p$ to $L^p$ for $p>1$:
$$\int\limits_{A_N}\sup_{(y,t)\in\Gamma_x}\left|f* P_t(y,t)\right|^{p_1}dx \leq C\int\limits_{A_N}(M |\tilde\alpha|^{\delta})^{\frac{p_1}{\delta}}+C\int\limits_{A_N}(M |\tilde\beta|^{\delta})^{\frac{p_1}{\delta}}\leq$$
$$\leq \int\limits_{\mathbb R^n} |\tilde \alpha|^{p_1}+C|A_N|^{1-\frac{p_1}{p_2}}\left(\int\limits_{A_N} \left(M|\tilde \beta|^{\delta}\right)^{\frac{p_2}{\delta}} \right)^{\frac{p_1}{p_2}}\leq C\alpha ^{p_1}+ C |A_N|^{1-\frac{p_1}{p_2}} \beta^{p_1}.$$
One can rewrite the line that we have just obtained in the following way:
$$C\alpha ^{p_1}+ C |A_N|^{1-\frac{p_1}{p_2}} \beta^{p_1}\ge \int\limits_{A_N}\sup_{(y,t)\in\Gamma_x}\left|f* P_t(y,t)\right|^{p_1}dx \ge \lambda^{p_1}|A_N|.$$
The last inequality is true, because $A_N \subset A.$
We deduce from here the following bound on $|A_N|$
$$\left|A_N\right|^{\frac{p_1}{p_2}}\leq C\left(\frac{\beta^{p_1}}{\lambda^{p_1}} + \left(\frac{\alpha^{p_1}}{\lambda^{p_1}}\right)^{\frac{p_1}{p_2}} \right).$$
After that the proof of the inequality for $|A|$ will be over just by tending $N$ to $\infty$.
\begin{lem}
Let $t>0$, $n>1$ and $0\leq t^n a+ tb-c$, where $a,b,c>0$ Then 
$$\frac{1}{t}\leq 2\left(\frac{b}{c}+ \left(\frac{a}{c}\right)^{\frac{1}{n}}\right).$$ 
\end{lem}
\begin{proof}
One can see that either $$t^n \frac{a}{c} \ge \frac{1}{2}$$ or $$t\frac{b}{c}\ge \frac{1}{2}.$$ From here we conclude that either $$\frac{1}{t}\leq 2 \frac{b}{c}$$ or $$\frac{1}{t}\leq \left(2\frac{a}{c}\right)^{\frac{1}{n}}.$$
In both cases
$$\frac{1}{t}\leq 2\left(\frac{b}{c}+ \left(\frac{a}{c}\right)^{\frac{1}{n}}\right).$$  
\end{proof}

Now, using that lemma with $t=\frac{1}{\left|A\right|^{\frac{p_1}{p_2}}} $,  $ c=\lambda^{p_1} $, $ a=\alpha^{p_1} $, $ b=\beta^{p_1}$, we will get the following inequality:
$$|A_N|^{\frac{p_1}{p_2}}\leq C\left(\frac{\beta^{p_1}}{\lambda^{p_1}}+\left(\frac{\alpha^{p_1}}{\lambda^{p_1}}\right)^{\frac{p_1}{p_2}}\right).$$
It remains to mention only that the needed estimate of $|A_N|$ follows from the above formula after writing $\lambda=\left(\frac{\beta^{p_2}}{\alpha^{p_1}}\right)^{\frac{1}{p_2-p_1}}$ and a little computation.

\end{proof}

Finally, after having estimated the measure $|A|$, we can perform the needed decomposition of $f$. First, observe that the following formulas hold true:
$$f(x)=\int\limits_{\mathbb R^{n+1}_+} \frac{\partial u}{\partial t}(y,t)\psi_t(x-y)dtdy,$$
where  $$u(y,t)=f*P_t(y,t), \psi_t=\frac{1}{t}\psi(\frac{x}{t}),$$ and $\psi$ is a smooth function such that $supp(\psi)\subset\{|x|<\frac{1}{2}\}.$
From here we deduce that
$$f(x)=\int\limits_{\bigcup\limits_{i\ge 1}\hat{A_i}}\frac{\partial u}{\partial t}(y,t)\psi_t(x-y)dtdy +\int\limits_{\mathbb R^n \backslash \bigcup\limits_{i\ge 1}\hat{A_i}}\frac{\partial u}{\partial t}(y,t)\psi_t(x-y)dtdy,$$
where $\hat A_i$ is the tent above $A_i.$
At last we are ready to define our functions $w$ and $v$.
$$w(x)=\int\limits_{\bigcup\limits_{i\ge 1}\hat{A_i}}\frac{\partial u}{\partial t}(y,t)\psi_t(x-y)dtdy,$$
$$v(x)=\int\limits_{\mathbb R^n \backslash \bigcup\limits_{i\ge 1}\hat{A_i}}\frac{\partial u}{\partial t}(y,t)\psi_t(x-y)dtdy).$$

Let us now prove the required properties of $w$ and $v$.

1).The estimate of the norm $\|v\|_{L^{p_2}}.$ 

Let $\phi \in L^{p_2\prime}(\mathbb R^n)$, $\|\phi\|=1,$
$$<v,\phi>=\int\limits_{\mathbb R^n}v\phi=\int\limits_{\mathbb R^n}\int\limits_{\mathbb R^n_+\backslash \hat A}\phi(x)\frac{\partial u}{\partial t}(y,t)\psi_t(x-y)dtdydx=\cdots,$$
where $\hat A:=\bigcup_{i\ge 1} \hat A_i$. Here, one can change the variables and use Holder's inequality:
$$\dots = \int\limits_{\mathbb R^n_+\backslash \hat A}\left(\int\limits_{\mathbb R^n}\phi(x)\psi_t(x-y)dx\right) \frac{\partial u}{\partial t}(y,t)dtdy = \int\limits_{\mathbb R^n_+\backslash \hat A}\phi*\psi_t(y)\frac{\partial u}{\partial t}(y,t)dtdy \leq$$
$$\leq  \left(\int\limits_{\mathbb R^{n+1}}\left|\phi*\psi_t(y)\right|^2 \frac{dy dt}{t}\right)^{\frac{1}{2}}\left(\int\limits_{\mathbb R^{n+1}\backslash \hat A} t \left| \nabla u \right|^2 dy dt\right)^{\frac{1}{2}}.$$
To estimate the first of the above integrals we use Parseval's identity:
$$\int\limits_{\mathbb R^{n+1}}\left|\phi*\psi_t(y)\right|^2 \frac{dy dt}{t}=\int\limits_{0}^{+\infty}\int\limits_{\mathbb R^{n}}\left|\phi*\psi_t(y)\right|^2 \frac{dy dt}{t}=\int\limits_{0}^{+\infty}\int\limits_{\mathbb R^{n}}\left|\hat\phi(\xi)\right|^2 \left|\hat\psi(t\xi)\right|^2 \frac{dy dt}{t}\leq$$
$$\leq C\|\phi\|_2^2.$$
For the second integral we use that $\Delta u = 0$ (because $u$ is harmonic):
$$\left|\nabla u\right|^2=\sum\limits_{i}\left(\frac{\partial u}{\partial x_i} \right)^2 =\frac{1}{2} \cdot 2 \sum\limits_{i}\left(\frac{\partial u}{\partial x_i} \right)^2 + \frac{1}{2} \cdot 2 \sum\limits_{i}\frac{\partial^2 u}{\partial x_i^2}=$$
$$ =\frac{1}{2} \sum\limits_{i}\frac{\partial^2 u^2}{\partial x_i^2}=\frac{1}{2} \Delta(u^2).$$
Therefore, using the second Green formula we get:
$$\int\limits_{\mathbb R^{n+1}_+\backslash \hat A} t \left| \nabla u \right|^2 dy dt=\frac{1}{2}\int\limits_{\mathbb R^{n+1}_+\backslash \hat A} t \Delta (u)^2 dy dt=$$
$$= \frac{1}{2}\int\limits_{\partial(\mathbb R^{n+1}_+\backslash \hat A)} \frac{\partial}{\partial n}(t) u^2 - \frac{\partial}{\partial n}(u^2) t dy dt \leq \dots ,$$
where $\frac{\partial}{\partial n}$ denotes the normal derivative.
$$\dots\leq \int\limits_{\partial(\mathbb R^{n+1}_+\backslash \hat A)}\left(t|u|\left|\frac{\partial u}{\partial n}\right| + \frac{1}{2}u^2\left|\frac{\partial t}{\partial n}\right|\right)=:\int\limits_{I}+\int\limits_{II}+\int\limits_{III}\left(t|u|\left|\frac{\partial u}{\partial n}\right| + \frac{1}{2}u^2\left|\frac{\partial t}{\partial n}\right|\right),$$
where 

$III:=\{(y,t)\in \mathbb R^{n+1}_+: |(y,t)|\leq R\}$ is a half-sphere in $\mathbb R^{n+1}_+$ and $R$ is some big number, 

$II:=\partial\left(\bigcup\limits_{i}\hat A_i\backslash\{(y,t): t=0\}\right)$, 

$I:=\partial\left(\bigcup\limits_{i}\mathbb R^{n+1}_+\backslash\hat A\right)\cap\{(y,t): t=0\}\backslash A.$

Now we are going to estimate consequently these three integrals. The first one can be easily bounded using the choice of $\lambda$

$$\int\limits_{I}\left(t|u|\left|\frac{\partial u}{\partial n}\right| + \frac{1}{2}u^2\left|\frac{\partial t}{\partial n}\right|\right)=\int\limits_{I}\frac{1}{2}u^2\left|\frac{\partial t}{\partial n}\right| \leq \int\limits_{I} \frac{1}{2}u^2\leq$$
$$\leq \int\limits_{\mathbb R^n\backslash A} \frac{1}{2}|f|^2\leq C\lambda^{2-p_1}\alpha^{p_1}=C \beta^2$$

The second one tends to zero as $R$ goes to infinity:

$$\int\limits_{III}\left(t|u|\left|\frac{\partial u}{\partial n}\right| + \frac{1}{2}u^2\left|\frac{\partial t}{\partial n}\right|\right) \rightarrow 0.$$

To establish the needed estimate of the second integral we first notice that $|u(x)|\leq \lambda$ for $x\in II.$ That is so, because each point of the corresponding tent is approximately at the same distance from $A$. Being more specific, we note that

$$\int\limits_{II}\left(t|u|\left|\frac{\partial u}{\partial n}\right| + \frac{1}{2}u^2\left|\frac{\partial t}{\partial n}\right|\right)\leq \int\limits_{II}C \lambda^2+ \int\limits_{II}t|u|\left|\frac{\partial u}{\partial n}\right|.$$

Let $x\in \hat A_k.$ There exists $x_0\in \mathbb R^n \backslash \bigcup\limits_{i\ge 1} A_i$ such that $dist(x_0,A_k)=|A_k|.$ That means that $x$ lies in the cone of the angle with the vertex at the point $x_0.$ From here we have $|u(x)|\leq \lambda.$

To end up the proof of the estimate, we shall now show that $t \left|\frac{\partial u}{\partial n}\right| \leq C\lambda$ for $(t,x)\in II.$  It is enough to show that $ t \left|\nabla u\right| \leq C\lambda$ for $(t,x)\in II.$
We provide the reader only with the estimate of the term $ t \left|\frac{\partial u}{\partial t}\right|,$ since other derivatives can be estimated the similar way.

Let us notice that there exists a ball with the center at $(t,x)$ and of radius $\frac{t}{C}$ that is contained in a cone with the vertex at $x_0$ for $x_0\in \mathbb R^n \backslash A.$
That is true, because from $(t,x)\in II$ we see that $t\ge \frac{edge(A_k)}{5}.$
Now, it is obvious that if we take a point $x_0\in \mathbb R^n \backslash A$ such that $dist(x,A_k)=5 diam(A_k)$ and a cone of an the aperture $30n$, then it contains the ball of radius $\frac{t}{2}$, centered at the point $(t,x).$ 
Armed with the existence of the mentioned ball, we are now able to use the following version of Green's Theorem:
$$\int\limits_{B_{x,t}}\frac{\partial u}{\partial t}= -\int\limits_{\partial B_{x,t}}u dx.$$
So, we deduce that
$$\left| \frac{\partial u}{\partial t}\right|=\left| \frac{2^n}{t^{n-1}}\int\limits_{B_{x,t}}\frac{\partial u}{\partial t}dx\right|=\left| \frac{2^n}{t^{n-1}}\int\limits_{\partial B_{x,t}} u dx\right|\leq C|u(x,t)|\leq C \lambda.$$ 
We are now finally ready to end up the estimate of our second integral:
$$\int\limits_{II}\left(t|u|\left|\frac{\partial u}{\partial n}\right| + \frac{1}{2}u^2\left|\frac{\partial t}{\partial n}\right|\right) \leq C \lambda^2 |A| +\lambda \int\limits_{II} t \left|\frac{\partial u}{\partial n}\right| \leq$$
$$\leq C\lambda^2 |A|+\lambda^2C|A|=C\lambda^2 |A|\leq C\beta^2.$$
As a consequence we get the following inequality
$$\int\limits_{\mathbb R^{n+1}_+\backslash \hat A} t \left| \nabla u \right|^2 dy dt\leq C\beta^2.$$
Finally, we obtain the needed bound of the scalar product $<v,\phi>$
$$<v,\phi>\leq C\|\phi\|_2\beta,$$
from where we get that
$$\|v\|_2\leq C\beta.$$

2). Now we need to deal with the function $w.$

First, it is obvious that the support of $w$ is the set $A.$

Second, let us bound $\int\limits_{\mathbb R^n}|w|^{p_1},$ using Holder's inequality and a trivial power mean  inequality:
$$\int\limits_{A}|w|^{p_1}=\int\limits_{A}|f-v|^{p_1}\leq C\int\limits_{\mathbb R^n}|f|^{p_1} + C\int\limits_{A}|v|^{p_1}\leq C\alpha^{p_1} + C|A|^{1-\frac{p_1}{p_2}}\beta^{p_2}\leq C\alpha^{p_1}.$$
We now estimate $\|w\|_{H^{p_1}}.$ To this end, we shall decompose $w$ into a sum of easier functions.
Recall that
$$w(x)=\int\limits_{\bigcup\limits_{i\ge 1}\hat A_i}\frac{\partial u}{\partial t}(y,t) \psi_t(x-y)dtdy,$$
where $A=\bigcup\limits_{A_i}$ is a Whitney decomposition of A, $\hat A_i$ is the tent  above $A_i.$

Consider the set $A^k$ defined the following way:
$$ A^k:=\{x: \sup\limits_{(y,t)\in \Gamma_x}|f*P_t(y,t)|>2^k \lambda\},$$
for $k=1,2,3,\dots$ and where we set $A_0:=A.$
Next, for each $A_k$ we consider it's Whitney decomposition:
$$A^k_i:=\bigcup\limits_{i\ge 1}A^k_i.$$
For each $k$ we take the union of sets of $A^k_i$ and pose:
$$\hat A^k:=\bigcup\limits_{i\ge 1}\hat A^k_i;$$$$ T_j^k:=\hat A^k_j \backslash \hat A^{k+1}.$$
It is obvious that 
$$\bigcup\limits_{j,k}T_j^k=\hat A^0=\hat A(=:\bigcup\limits_{i\ge 1}\hat A_i).$$

We are able now to write out the corresponding decomposition:
$$w(x)=\int\limits_{\bigcup\limits_{i\ge 1} \hat A_i}\frac{\partial u}{\partial t}(y,t) \psi_t(x-y) dt dy=\sum\limits_{j=1}^{\infty}\sum\limits_{k=0}^{\infty}\int\limits_{T_j^k}\frac{\partial u}{\partial t}(y,t) \psi_t(x-y) dt dy.$$
We set:
$$g_{j,k}(x)=\int\limits_{T_j^k}\frac{\partial u}{\partial t}(y,t) \psi_t(x-y) dt dy.$$
It is obvious that the support of $g_{j,k}$ is in $A_j^k.$ One can also see that holds $\int\limits_{T_j^k} g_{j,k}(x) dx=0.$

Let us now prove the bound on $g_j^k$, that we will need to end up the estimate of $w$, i.e.
$$\|g_j^k\|_2\leq C|A_{j,k}|^{\frac{1}{2}}2^k \lambda.$$
We act much the same way as while getting the bound on the function $v.$

Consider $\phi \in L^2(\mathbb R^n)$, such that $\|\phi\|_2=1$. We will now estimate $<g_{j,k},\phi>$. 
$$<g_{j,k},\phi>=\int\limits_{\mathbb R^n}g_{j,k}\phi\leq \left(\int\limits_{\mathbb R^{n+1}}\left|\phi*\psi_t(y)\right|^2 \frac{dy dt}{t}\right)^{\frac{1}{2}}\left(\int\limits_{T_{j,k}} t \left| \nabla u \right|^2 dy dt\right)^{\frac{1}{2}}\leq$$
$$C\|\phi\|_2\left(\int\limits_{T_{j,k}} t \left| \nabla u \right|^2 dy dt\right)^{\frac{1}{2}}\leq C||\phi||_2 \left(\int\limits_{\partial T_{j,k}}t|u|\left|\frac{\partial u}{\partial n}\right| + \frac{1}{2}u^2\left|\frac{\partial t}{\partial n}\right|\right)^{\frac{1}{2}}\leq \dots$$
The estimates above are much the same as in the bound of $<v,\phi>.$ We leave them without any further comments.

Notice that $\left|\frac{\partial t}{\partial n}(x)\right|\leq C$ for $x\in \partial T_{j,k}.$ It can be proved using the same methods as in the estimates of $v$ that $|u(x)|\leq 2^{k+1}\lambda$ for $x\in \partial T_{j,k}$ and that $t\frac{\partial u}{\partial n}\leq 2^{k+1} \lambda C$ for $x\in \partial T_{j,k}$.
From the written above, we conclude that:
$$C\|\phi\|_2 \left(\int\limits_{\partial T_{j,k}}t|u|\left|\frac{\partial u}{\partial n}\right| + \frac{1}{2}u^2\left|\frac{\partial t}{\partial n}\right|\right)^{\frac{1}{2}}\leq C(|\partial T_{j,k}| 2^{2k}\lambda^2)^{\frac{1}{2}}\leq C|A_{j,k}^{\frac{1}{2}}|2^k \lambda.$$

Finally, we end up the estimate of $<g_{j,k},\phi>:$
$$\dots \leq C\|\phi\|_2|A_{j,k}|^{\frac{1}{2}}2^k \lambda.$$
As a result we have:
$$w(x)=\sum\limits_{j,k}g_{j,k}(x)=\sum\limits_{j,k}\lambda_j^k a_{j,k}(x),$$
where $\lambda_j^k:=C2^k\lambda |A_{j,k}|^{1/p_1}$ and $a_{j,k}:=\frac{g_j^k}{\lambda_j^k}$.
From here and from the properties of $g_{j,k}$ we find that $a_{j,k}$ are in fact Hardy atoms. 
$$\|w\|^{p_1}_{H^{p_1}(\mathbb R^n)}\leq C \sum\limits_{j,k}|\lambda_j^k|^{p_1}=C\sum\limits_{j,k}2^{kp_1}\lambda^{p_1} |A_{j,k}|=$$
$$=C\lambda^{p_1}\sum\limits_{k=0}^{\infty}2^{kp_1}(\sum\limits_{j}|A_{j,k}|)=C\lambda^{p_1}\sum\limits_{k=0}^{\infty}2^{kp_1}|A^k|\leq C \int\limits_{A_0}\sup\limits_{(y,t)\in \Gamma_x}|f*P_t(y,t)|^{p_1}dx\leq \dots$$

From the fact that 
$$A^k=\{x : \sup\limits_{(y,t)\in \Gamma_x}|f*P_t(y,t)|>2^k \lambda\}$$
we deduce the following identity
$$|A^k|=|{x : \sup\limits_{(y,t)\in \Gamma_x}|f*P_t(y,t)|^{p_1}>2^{k p_1} \lambda^{p_1}}|.$$
Observing the inequalities
$$2^{(k+1)p_1}\lambda^{p_1} |A^{k+1}|\leq \frac{2^{p_1}}{2^{p_1}-1}\left(2^{(k+1)p_1}-2^{kp_1}\right)|A^{k+1}|$$
one can easily sum them up and obtain
$$\lambda^{p_1}\sum\limits_{k=0}^{\infty}2^{kp_1}|A^k|\leq C \int\limits_{A_0}\sup\limits_{(y,t)\in \Gamma_x}|f*P_t(y,t)|^{p_1}dx.$$

At last, we end up the estimate:
$$\dots = C \int\limits_{A}\sup\limits_{(y,t)\in \Gamma_x}|f*P_t(y,t)|^{p_1}dx\leq C\alpha^{p_1}+C|A|^{1-\frac{p_1}{p_2}}\beta^{p_1}\leq C\alpha^{p_1}.$$
\end{proof}
\begin{rem}
One can mention that the proof of the Theorem 1 is closely related to atomic decomposition of functions in $H^p$ spaces for $p\leq1$ (see \cite{GCRF} for the details).
\end{rem}
\begin{rem}
It is natural to conjecture that the bound $\frac{n-1}{n}$ is not sharp. Alas, in the case of the spaces $H^p(\mathbb R^n)$ where $p<\frac{n-1}{n}$ the subharmonicity argument applied above breaks down. This is so because the system of partial differential equations that defines those Hardy space is way more complicated  (see \cite{KG}).
\end{rem}
\renewcommand{\refname}{References}

\end{document}